\documentclass[a4paper,11pt,reqno,noindent]{amsart}
\usepackage[centertags]{amsmath}
\usepackage{amsfonts,amssymb,amsthm} 
\usepackage[english]{babel}
\usepackage{newlfont}
\usepackage{color}
\usepackage[body={15cm,21.5cm},centering]{geometry} 
\usepackage{fancyhdr}
\usepackage{esint}

\usepackage{enumerate}

\newtheorem{theo}{Theorem}
\newtheorem{lemma}{Lemma}

\newtheorem{coro}{Corollary}

\theoremstyle{definition}
\newtheorem{rem}{Remark}
\newtheorem{defi}{Definition}

\numberwithin{equation}{section}
\numberwithin{lemma}{section} 
\numberwithin{defi}{section} 
\numberwithin{rem}{section} 

\newcommand \dps{\displaystyle }

\newcommand{\ho}{\mathrm{hom}}

\newcommand{\R}{\mathbb{R}}
\newcommand{\Z}{\mathbb{Z}}
\newcommand{\N}{\mathbb{N}}
\newcommand{\T}{\mathbb{T}}

\newcommand{\e}{\varepsilon}

\newcommand{\ee}{\mathbf{e}}

\mathchardef\emptyset="001F

\def\aa{\boldsymbol a}
\newcommand{\bb}{\boldsymbol b}

\newcommand{\dig}[1]{\mathrm{diag}\left[ #1\right]}

\newcommand{\cov}[2]{\mathrm{cov}\left[#1;#2\right]}
\newcommand{\expec}[1]{\left\langle #1 \right\rangle}
\newcommand{\Expec}[1]{\left\langle #1 \right\rangle}
\newcommand{\step}[1]{\noindent \textit{Step} #1.}

\newcommand{\osc}[2]{\underset{\dps #1}{\mathrm{osc}} \,#2\,}

\title[Quantitative expansion in stochastic homogenization]
{An optimal quantitative two-scale expansion in stochastic homogenization of discrete elliptic equations}
\author[A. Gloria, S. Neukamm \& F. Otto]{Antoine Gloria, Stefan Neukamm \& Felix Otto}
\date{\today}
\address[Antoine Gloria]{Universit\'e Libre de Bruxelles (ULB) \\ Brussels, Belgium \\ and Project-team SIMPAF \\  Inria Lille - Nord Europe \\ Villeneuve d'Ascq, France}
\email{agloria@ulb.ac.be}
\address[Stefan Neukamm]{Max-Planck-Institut f\"ur Mathematik in den Naturwissenschaften \\ Leipzig, Germany}
\email{neukamm@mis.mpg.de}
\address[Felix Otto]{Max-Planck-Institut f\"ur Mathematik in den Naturwissenschaften \\ Leipzig, Germany}
\email{otto@mis.mpg.de}
\begin{document}
\maketitle

\begin{center}
\begin{minipage}{13cm}
\small{
\noindent {\bf Abstract.} 
We establish an optimal, linear rate of convergence for the stochastic homogenization of discrete linear elliptic equations. We consider the model problem of independent and identically distributed coefficients on a discretized unit torus. We show that the difference between the solution to the random problem on the discretized torus and the first two terms of the two-scale asymptotic expansion has the same scaling as in the periodic case. In particular the $L^2$-norm in probability of the \mbox{$H^1$-norm} in space of this error scales like $\e$, where $\e$ is the  discretization parameter of the unit torus.
The proof makes extensive use of previous results by the authors, and  of recent annealed estimates on the Green's function by Marahrens and the third author.
\vspace{10pt}


\vspace{6pt}
}

\end{minipage}
\end{center}

\bigskip


\section{Introduction}

\noindent We establish a linear rate of convergence for the stochastic homogenization of discrete linear elliptic equations, which is optimal.
Before we turn to the stochastic case, let us recall some standard analogous results in the periodic case.
Let $A$ be a uniformly elliptic and bounded symmetric matrix field on
the unit torus $\T:=(\R/\Z)^d$, $D\subset\R^d$ be a smooth domain, and $f$ be a smooth function.
Let $\e>0$. It is well-known that the unique weak solution $u_\e \in H^1_0(D)$ of the linear elliptic equation
\begin{equation}\label{eq:I:1}
\left\{
\begin{array}{rcl}
-\nabla \cdot A\Big(\frac{\cdot}{\e}\Big)\nabla u_\e &=&f \mbox{ in }D,\\
u_\e&=&0 \mbox{ on }\partial D
\end{array}
\right.
\end{equation}
converges weakly in $H^1(D)$ as $\e \to 0$ to the unique weak solution $u_{\hom}\in H^1_0(D)$ of the homogenized equation
\begin{equation}
  \left\{
    \begin{array}{rcl}
      -\nabla \cdot A_\ho\nabla u_{\hom} &=&f \mbox{ in }D,\\
      u_{\hom}&=&0 \mbox{ on }\partial D.
    \end{array}
  \right.
\end{equation}
The homogenized matrix $A_\ho$ is symmetric and characterized for all $\xi \in \R^d$ by
\begin{equation*}
  \xi\cdot A_\ho \xi\,=\, \int_{\T} (\nabla \phi_\xi+\xi)\cdot A (\nabla \phi_\xi+\xi),
\end{equation*}
where $\phi_\xi$ is the unique  weak solution in $H^1(\T)$ with vanishing mean of the periodic corrector equation in direction $\xi$:
\begin{equation*}
-\nabla \cdot A (\nabla \phi_\xi+\xi)\,=\,0 \mbox{ in }\T.
\end{equation*}
From a two-scale expansion, one formally expects for all $x$ in the interior of $D$:
\begin{equation}\label{intro:eq:asexp}
  u_\e(x)\,=\,u_{\hom}(x)+\e \sum_{j=1}^d\phi_{j}(\tfrac{x}{\e})\nabla_j u_{\hom}(x) +o(\e),
\end{equation}
where the $\phi_j$ are the correctors in the canonical directions $\ee_j$ (extended by periodicity to $\R^d$).
This identity cannot hold at the boundary since the correctors
$\phi_j$ do not satisfy the homogeneous Dirichlet boundary conditions:
There is indeed a boundary layer.
Yet, for all $\tilde{D}$ compactly supported in $D$, Allaire and Amar \cite[Theorem~2.3]{Allaire-Amar-99} proved the following rigorous version of \eqref{intro:eq:asexp} if the coefficients $A$ are H\"older continuous:
\begin{equation}\label{intro:per-H1}
  \|u_\e-u_{\hom}-\e \sum_{j=1}^d\phi_j(\tfrac{\cdot}{\e})\nabla_j u_{\hom} \|_{H^1(\tilde{D})} \,\leq \, C\e ,
\end{equation}
where the multiplicative constant $C$ only depends on
$\tilde{D}$, $D$, the H\"older exponent and norm of $A$, and on the $C^3$-norm of $u_{\hom}$ (which is smooth since $f$ and $D$ are smooth).
This result relies on a previous work by Avellaneda and Lin \cite[Theorem~5]{Avellaneda-Lin-87}, who proved under similar assumptions that
\begin{equation}\label{intro:per-L2}
\|u_\e-u_{\hom}\|_{L^\infty ({D})} \,\leq \, C\e.
\end{equation}
These are bounds on the homogenization error.
\smallskip

\noindent
In the stochastic case, we consider a matrix field $A$ that is stationary and ergodic, in place of periodic. We refer the reader to \cite{Papanicolaou-Varadhan-79} for details. In order to obtain quantitative estimates in the spirit of \eqref{intro:per-H1} and \eqref{intro:per-L2}, one has to make assumptions on the statistics of $A$ in addition to ergodicity.
There are few results in the literature on quantitative estimates of the homogenization error for elliptic equations in divergence form in the stochastic case.
In \cite[Theorem~3.1]{Yurinskii-86}, Yurinski{\u\i} proved for algebraically decaying correlations that for all $d>2$, there exists some H\"older exponent $\gamma>0$ and a function $T$ of $\e$ such that
\begin{equation}\label{intro:Yu}
\expec{\|u_\e-u_{\hom}-\e \sum_{j=1}^d\phi_{T(\e),j}(\tfrac{\cdot}{\e})\nabla_j u_{\hom}\|^2_{H^1(D)}}^{1/2}\,\leq\, C\e^\gamma,
\end{equation}
where $\phi_{T,j}$ is the modified corrector, which is the stationary, almost sure solution to 
\begin{equation}
  T^{-1} \phi_{T,j}-\nabla \cdot A(\nabla \phi_{T,j}+\ee_j)\,=\,0 \mbox{ in }\R^d.
\end{equation}
This equation is an approximation of the corrector equation when $T\to \infty$.
This is the first quantitative result in stochastic homogenization.
Note that a formal linearization in the case of small ellipticity contrast $\lambda \uparrow 1$ yields $\gamma=1$ for $d>2$.
Besides not covering dimension $d=2$, the work by Yurinski{\u\i} does not allow to reach the scaling $\gamma=1$, even in the case of small ellipticity constrast  (and for domains $\tilde D$ compactly included in $D$).

\medskip

\noindent
In this article, we  simplify the above setting with three respects:
\begin{enumerate}[(a)]
\item We treat \textbf{finite-difference equations} instead of continuum partial differential equations. More precisely, we replace
  the differential operator in \eqref{eq:I:1} by the finite difference operator
  $\nabla^*_\e\aa(\tfrac{\cdot}{\e})\nabla_\e$ which acts on discrete
  functions defined on the scaled lattice $\e\Z^d$, $\e>0$. Here $\nabla_\e$ and $\nabla^*_\e$ are finite difference
  approximations of the continuum gradient $\nabla$ and continuum
  divergence operator $-\nabla\cdot $, and are defined for scalar fields $v:\e\Z^d\to\R$ and vector fields, $g=(g_1,\ldots,g_d):\e\Z^d\to\R^d$ as follows:
  \begin{align}
    \label{eq:nablaeps}
    \begin{split}
      &\nabla_{\e,i} v(y) = \e^{-1}(v(y+\e\ee_i)-v(y)), \quad
      \nabla^*_{\e,i} v(y) =
      \e^{-1}(v(y-\e\ee_i)-v(y)),\\
      &\nabla_\e v=(\nabla_{\e,1}v,\ldots,\nabla_{\e,d}v),\qquad
      \nabla^*_\e g=\sum_{i=1}^d\nabla^*_{\e,i}g_i.
    \end{split}
  \end{align}
  The coefficients of the operator $\nabla^*_\e\aa(\tfrac{\cdot}{\e})\nabla_\e$ are given by a
  randomly chosen matrix field \mbox{$\aa:\Z^d\to\R^{d\times d}$}. We
  assume that $\aa$ is diagonal and uniformly elliptic in the sense
  that 
  \begin{equation*}
    \aa(x)\in\Omega_0:=\left\{\,\left(
      \begin{array}{ccc}
        \aa_1&&0\\&\ddots&\\0&&\aa_d
      \end{array}
\right)\,\Big|\,\lambda\leq
    \aa_j\leq 1,\;j=1,\ldots,d\,\right\}
  \end{equation*}
  for all $x\in\Z^d$. Above, $\lambda>0$ is a (deterministic)
  ellipticity constant  and fixed throughout the paper.
\item In order to avoid boundary layers, we consider a problem on the
  (discretized) unit \textbf{torus} $\T_\e:=(\e\Z/\Z)^d$ and assume
  w.~l.~o.~g. that the grid size $\e$ is the inverse of a large integer $L:=\e^{-1}\in \N$.
  The equation under consideration is the following finite difference equation with random
  coefficients
  \begin{equation}\label{eq:torus-L}
    \nabla_\e^*\aa(\tfrac{\cdot}{\e})\nabla_\e u_\e\,=\,f_\e \mbox{ on }\T_\e,\qquad\sum_{\T_\e}u_\e=0,
  \end{equation}
  where $f_\e\,:\,\T_\e\to\R$ is a deterministic r.~h.~s. with zero
  mean, and should be viewed as a discretization of a continuum
  r.~h.~s. $f\in L^2(\T)$. Note that on $\T_\e$ we have by periodicity the \textit{discrete integration by parts formula}
  \begin{equation*}
    \sum_{\T_\e}\nabla_\e v\cdot g=\sum_{\T_\e}v \nabla^*_\e g\qquad\text{for all $v:\T_\e\to\R$ and $g:\T_\e\to \R^d$}.
  \end{equation*}

\item We consider the simplest possible statistics, namely \textbf{independent
    and identically distributed} (i.~i.~d.) coefficients. Since \eqref{eq:torus-L} is posed on the discretized unit torus $\T_\e$, $\e^{-1}=L\in \N$, we in fact consider the \textit{periodic
  i.\ i.\ d.\ ensemble} which is constructed as follows: For a fixed measure $\beta$ on our matrix space $\Omega_0$ (``single-site measure'') and for
  fixed size $L$ of the discrete, rescaled torus $\mathbb{T}_L:=(\Z/L\Z)^d=\frac{1}{\e}\T_\e$, 
  we shall throughout the paper denote by $\langle\cdot\rangle$ the \textit{periodic
  i.\ i.\ d.\ ensemble} on 
  $$\Omega_L=\{\aa\colon\mathbb{Z}^d\rightarrow\Omega_0\,|\,\forall z\in\mathbb{Z}^d\;\aa(\cdot+Lz)=\aa\}
  \,\widehat = \,\Omega_0^{\mathbb{T}_L}.$$ For a random variable, i.\ e.\ a measurable function
  $\zeta\colon\Omega_L\rightarrow\mathbb{R}$, it is given by the product measure
  \begin{equation}
    \expec{\zeta}=\prod_{x\in\mathbb{T}_L}\int_{\Omega_0}\zeta(\aa)\,\beta(d\aa(x)).\label{eq:def-iidens}
  \end{equation}
  Evidently, our i.\ i.\ d.\ ensemble is {\it stationary} in the sense that for all shifts 
  $y\in\mathbb{Z}^d$ the random (periodic) tensor fields $\aa$ and $\aa(\cdot+y)$ 
  have the same distribution in $\Omega_L$.
\end{enumerate}
\noindent
In the above setting the qualitative theory of stochastic
homogenization applies (see e.~g. Kozlov~\cite{Kozlov-79}, Papanicolaou and Varadhan
\cite{Papanicolaou-Varadhan-79} and K\"unnemann \cite[Theorem~4]{Kunnemann-83}): There exists a deterministic, symmetric, positive
definite $d\times d$ matrix $\aa_{\ho}$ (only depending on the
single-site probability measure $\beta$) such
that the following statement is true: Suppose that $f_\e$ converges (in a discrete $H^{-1}$-norm) to some
function $f\in L^2(\T)$ with $\int_{\T}f=0$. Let $u_\e$ be the unique
(random) solution with vanishing mean to
\eqref{eq:torus-L}, and let $u_{\hom}\in H^1(\T)$ be the unique
solution to the continuum homogenized equation
\begin{equation}\label{eq:hom-eq-torus}
  -\nabla\cdot \aa_\ho \nabla u_{\hom} \,=\,f \mbox{ in }\T,\qquad\int_{\T}u_{\hom}=0.
\end{equation}
Then $u_\e$ converges almost surely (weakly in the discrete $H^1$-norm) to $u_{\hom}$.
\medskip

\noindent
We are interested in estimates on the speed of convergence of $u_\e$
to $u_{\hom}$. The error $u_\e-u_{\hom}$ consists of two contributions that are of different nature, namely the
\textit{homogenization error} and the \textit{discretization error}. The latter is purely deterministic and introduced by approximating the deterministic, \emph{continuum}, homogenized equation \eqref{eq:hom-eq-torus} by the deterministic, homogenized, \emph{finite difference} equation
\begin{equation}\label{eq:def-fL}
  \nabla^*_\e\aa_\ho \nabla_\e u_{0,\e} \,=\,f_\e \mbox{ in }\T_\e,\qquad\sum_{\T_\e}u_{0,\e}=0,
\end{equation}
where $f_\e\in L^2(\T_\e)$ denotes a suitable approximation of the
continuum r.~h.~s. $f$. The discretization error is
well-studied and well-understood (see for instance \cite{Leveque-07}).

\medskip

\noindent
In this paper we focus on the homogenization error that monitors the
difference between the random, variable-coefficient equation \eqref{eq:torus-L} and the deterministic, constant-coefficient equation \eqref{eq:def-fL}.
The main result of this paper is the upcoming quantitative two-scale
expansion estimate. We quantify the error
in the discrete $L^2(\T_\e)$- and $H^1(\T_\e)$-norms, which are defined for $v:\T_\e \to \R$ as
\begin{equation*}
  \|v\|_{L^2(\T_\e)}=\sqrt{\e^d\sum_{\T_\e}v^2},\qquad \|v\|_{H^1(\T_\e)}=\sqrt{\e^d\sum_{\T_\e} \left(v^2+|\nabla_\e v|^2\right)}.
\end{equation*}
\medskip

\begin{theo}\label{th:main}
  Let $d\geq 2$ and $a$ be i.~i.~d. coefficients.
Then there exists a deterministic, symmetric, positive definite matrix $\aa_{\hom}\in\R^{d\times d}$ (only
  depending on $\beta$ and $d$) with the following property.
  \smallskip
  
  \noindent
  Given $\e>0$ with $\e^{-1}=L\in\N$, and a r.~h.~s.
  \begin{equation*}
    f_\e:\T_\e\to\R\qquad\text{with }\sum_{\T_\e}f_\e=0,
  \end{equation*}
  let $u_\e:\Omega_L\times\T_\e\to\R$ and $u_{0,\e}:\T_\e\to\R$ be the
  unique solutions to \eqref{eq:torus-L} and
    \eqref{eq:def-fL}, respectively. Then
    \begin{multline}\label{eq:main}
      \expec{\|u_\e-u_{0,\e}-\e\sum_{j=1}^d\phi_j(\tfrac{\cdot}{\e})\nabla_{\e,j}
        u_{0,\e} \|^2_{H^1(\T_\e)}}^{1/2} \\\,\lesssim \,\e\, ||f_\e||_{L^2(\T_\e)}\,\left\{
        \begin{array}{ll}
          (\ln \frac{1}{\e})^{1/2} &\mbox{ for }d=2
          \\
          1 &\mbox{ for }d>2
        \end{array}
      \right.
      ,
    \end{multline}
    where $\phi_1,\ldots,\phi_d$ are the periodic correctors
    associated with the periodic i.~i.~d. ensemble $\expec{\cdot}$ via
    the periodic corrector equation (see \eqref{eq:per-corr} below). The
    multiplicative constant in \eqref{eq:main} only depends on the
    constant of ellipticity $\lambda$,
    and the dimension $d$.
\end{theo}
\noindent As a corollary we get:
\medskip

\begin{coro}\label{cor:1}In the situation of Theorem~\ref{th:main} we have
\begin{eqnarray}
  \expec{\|u_\e-u_{0,\e}\|^2_{L^2(\T_\e)}}^{1/2} &\lesssim & \e\,\|f_\e\|_{L^2(\T_\e)}\, \left\{
\begin{array}{ll}
(\ln \tfrac{1}{\e})^{1/2} &\mbox{ for }d=2
\\
1 &\mbox{ for }d>2
\end{array}
\right. ,\label{eq:coro1}
\end{eqnarray}
The multiplicative constant in this estimate only depends on $\lambda$ and $d$.
\end{coro}
\noindent Note that a bound similar to \eqref{eq:coro1} (with however a suboptimal exponent) was recentlty obtained by Conlon and Spencer  in \cite{Conlon-Spencer-12}, in the case of the whole space and a massive term. This result complements the estimates of the quantity $u_\e-\expec{u_\e}$ obtained in \cite[Theorem~2]{Gloria-12c}, and in their optimal form in \cite[Corollary~2 and Corollary~3]{Marahrens-Otto-13}. In particular, the annealed estimates of \cite[Theorem~1]{Marahrens-Otto-13} (see Lemma~\ref{L:MO} below) on the gradient of the (periodic) Green's function, combined with the argument leading to \cite[Theorem~2]{Gloria-12c}, allow one to prove that 
$$
\expec{\|u_\e-\expec{u_\e}\|_{L^2(\T_\e)}^2}^{1/2} \,\lesssim \, 
\e\,\|f_\e\|_{L^{2}(\T_\e)}\, \left\{
\begin{array}{ll}
(\ln \tfrac{1}{\e})^{1/2} &\mbox{ for }d=2
\\
1 &\mbox{ for }d>2
\end{array}
\right. .
$$
The combination of this estimate with Corollary~\ref{cor:1} then yields an estimate of the difference between the solution $u_{0,\e}$ of the problem with constant coefficients and the expectation $\expec{u_\e}$ of the solution of the original problem --- i.~e. the systematic error:
$$
\|u_{0,\e}-\expec{u_\e}\|_{L^2(\T_\e)} \,\lesssim \, 
\e\,\|f_\e\|_{L^{2}(\T_\e)}\, \left\{
\begin{array}{ll}
(\ln \tfrac{1}{\e})^{1/2} &\mbox{ for }d=2
\\
1 &\mbox{ for }d>2
\end{array}
\right. .
$$
It is rather surprising that we have to go through the $H^1$-norm in order to control this systematic error, or conversely that we do not have to estimate this term in order to prove the (seemingly stronger) statement of Theorem~\ref{th:main}.
Compared to the work \cite{Yurinskii-86} by Yurinski{\u\i}, Theorem~\ref{th:main} covers dimension $d=2$ and gives optimal estimates (in terms of scaling in $\e$) in any dimension, as can be seen by considering the regime of small ellipticity contrast. 
Note that for $d=1$, the scaling is different: \eqref{eq:coro1} is expected to be replaced by
\begin{equation*}
  \expec{\|u_\e-u_{0,\e}\|^2_{L^2(\T_\e)}}^{1/2} \,\lesssim \, \sqrt{\e}\,\|f_\e\|_{L^2(\T_\e)},
\end{equation*}
as explicity checked in the continuum setting by Bourgeat and Piatnitski in \cite{Bourgeat-99}.
As opposed to these works, the present analysis heavily relies on the use of a spectral gap estimate in the probability space. We refer the reader to \cite{Gloria-Neukamm-Otto-11a} for relevant references on the subject.

\begin{rem}
  The reason why we consider the discrete setting is the
  following: In the proof of Theorem~\ref{th:main} we make extensive use of
  recent, quantitative results that we obtained in a series of paper  \cite{Gloria-Otto-09,Gloria-Otto-09b,
    Gloria-Neukamm-Otto-11a,Marahrens-Otto-13} in the discrete setting. The extension of some results to the continuum setting is currently under investigation, see \cite{Gloria-Otto-10b}.
\end{rem}

\begin{rem}(Rescaling).\label{R:rescaling}
  For the 
proof of Theorem~\ref{th:main} it is
  convenient to rescale the discretized torus $\T_\e$ so that the
  grid size becomes $1$. Recall that
  $\e=\frac{1}{L}$ for some integer $L\in\N$. Hence, 
the $L$-rescaled version of $\T_\e$ yields the discrete rescaled torus
  $\T_L=(\Z/L\Z)^d=\frac{1}{\e}\T_\e$.
In analogy with \eqref{eq:nablaeps}  we introduce discrete derivatives $\nabla$ and $\nabla^*$
acting on scalar fields $v:\T_L\to\R$ and
  vector fields $g=(g_1,\ldots,g_d):\T_L\to\R^d$ as follows:
  \begin{align}\label{eq:def-disder}
    \begin{split}
      &\nabla_i v(x) = v(x+\ee_i)-v(x), \quad \nabla^*_i v(x) =
      v(x-\ee_i)-v(x),\\
      &\nabla v=(\nabla_1v,\ldots,\nabla_dv),\qquad
      \nabla^*g=\sum_{i=1}^d\nabla^*_ig_i.
    \end{split}
  \end{align}
  In order to state Theorem~\ref{th:main} in its rescaled
  version, we set
\begin{equation*}
  u(x):=u_\e(\e x),\qquad u_{0}(x):=u_{0,\e}(\e x),\qquad
  \tilde f(x):=\e^{2}f_\e(\e x)\qquad\text{for }x\in\Z^d.
\end{equation*}
So defined, $u$, $u_{0}$ and $\tilde f$ are functions on
the rescaled torus $\T_L$ with vanishing mean, and the finite difference equations \eqref{eq:torus-L} and \eqref{eq:def-fL} turn
into 
\begin{eqnarray}\label{eq:9}
  \nabla^*\aa\nabla u&=&\tilde f\mbox{ on }\T_L,\qquad \sum_{\T_L}u  =0,\\\label{eq:9b}
  \nabla^*\aa_{\hom}\nabla u_{0}&=&\tilde f\mbox{ on }\T_L,\qquad \sum_{\T_L}u_{0}=0.
\end{eqnarray}
Furthermore, the two-scale expansion in \eqref{eq:main} takes the form
\begin{equation}\label{eq:uzhat}
  z:=u  -u_{0}-\sum_{j=1}^d\phi_j\,\nabla_j
  u_{0},
\end{equation}
and estimate \eqref{eq:main} of Theorem~\ref{th:main} can be reformulated
as
\begin{equation}\label{eq:13}
  \expec{\sum_{\T_L}(z^2+L^2|\nabla
    z|^2)}^{\frac{1}{2}}\lesssim
  L\mu_d^{\frac{1}{2}}(L)\,\left(\sum_{\T_L}\tilde f^2\right)^{\frac{1}{2}},
\end{equation}
where we set for abbreviation
\begin{equation*}
  \mu_d(L)=\begin{cases}
    \ln L&\text{for }d=2,\\
    1&\text{for }d>2.
  \end{cases}
\end{equation*}
In fact we shall establish \eqref{eq:13} (and thus \eqref{eq:main}) in form of the estimate
\begin{equation}\label{eq:mainL}
  \expec{\sum_{\T_L}(z^2+L^2|\nabla
    z|^2)}^{\frac{1}{2}}\lesssim L\mu_d^{\frac{1}{2}}(L)\,\left(\sum_{\T_L}|\nabla^2u_0|^2\right)^{\frac{1}{2}},
\end{equation}
where $\nabla^2 u_0(x)$ denotes the discrete Hessian of $u_0$ at $x$
and is given by the $d\times d$ matrix with entries $-\nabla^*_i\nabla_j
u_0(x)$. Note that \eqref{eq:mainL} indeed implies \eqref{eq:13},
since $u_0$ (as a solution to the constant-coefficient
difference equations \eqref{eq:9b}) satisfies
the a priori estimate $\sum_{\T_L}|\nabla^2
u_0|^2\lesssim \sum_{\T_L}f^2$ up to a multiplicative constant that
only depends on $\lambda$ and $d$. 
\end{rem}
\medskip

\subsection*{Notation}
Throughout this article, we use the following notation:
\smallskip

\begin{itemize}
\item $d\geq 2$ is the dimension;
\item $(\ee_1,\dots,\ee_d)$ denotes the canonical basis of $\Z^d$; 
\item $\T_L=(\Z/L\Z)^d$ denotes the discretized  $L$-rescaled torus;
\item $x\mod L$ denotes the unique point in  $([0,L)\cap\Z)^d$ with
  $x=(x\mod L)+Lx'$ for some $x'\in \Z^d$;
\item $\lesssim$ and $\gtrsim$ stand for $\leq$ and $\geq$ up to a multiplicative constant which only depends on the quantities specified in the context;
\item when both $\lesssim$ and $\gtrsim$ hold, we simply write $\sim$;
\item $\beta$ denotes a single-site probability measure on
$\Omega_0$, see Section~\ref{sec:assumptions-ensemble};
\item $\expec{\cdot}$ denotes the $L$-periodic i.~i.~d. ensemble on $\Omega_L$ associated with $\beta$, see Section~\ref{sec:assumptions-ensemble};
\item $\cov{\cdot}{\cdot}$ denotes the covariance associated with $\expec{\cdot}$;
\item we denote the $(i,j)$-th entry of a $d\times d$-matrix $\bb$ by $\bb^{ij}$ and write ``$:$''  for the inner
  product in $\R^{d\times d}$, i.~e. \mbox{$\aa:\bb=\sum_{i,j=1}^d\aa^{ij}\bb^{ij}$};
\item for all $L>0$, $\mu_d(L)=\ln L$ for $d=2$ and $\mu_d(L)=1$ for $d>2$.
\end{itemize}

\section{Assumptions on the ensemble and the notion of the corrector}
\noindent
In this section we introduce and motivate the requi
 assumptions on the ensemble, the definition of the corrector and the homogenized coefficients. We recall some recent quantitative estimates from stochastic homogenization that are at the basis of the proof of Theorem~\ref{th:main}. Finally, we comment on the role played by the i.~i.~d. assumption.
\subsection{Assumptions on the ensemble}\label{sec:assumptions-ensemble}
\noindent
Recall that $\expec{\cdot}$ denotes the i.~i.~d. ensemble associated
with the single-site measure $\beta$ via \eqref{eq:def-iidens}. Our theory involves only two probabilistic ingredients: 
A covariance estimate and a Logarithmic Sobolev Inequality (LSI) for $\langle\cdot\rangle$. 
We start with the covariance estimate, which is explicitly used in the proof of this paper.
We shall comment on the LSI in the following section.

\begin{lemma}\label{L1} Let $\langle\cdot\rangle$ denote the periodic
  i.\ i.\ d.\ ensemble (see \eqref{eq:def-iidens}). Then
  we have for any two random variables $\zeta$, $\tilde\zeta$:
  \begin{equation}\label{O.1}
    \cov{\zeta}{\tilde\zeta}\le\sum_{y\in\mathbb{T}_L}
    \langle(\frac{\partial      \zeta}{\partial y})^2\rangle^\frac{1}{2}
    \langle(\frac{\partial\tilde\zeta}{\partial y})^2\rangle^\frac{1}{2}.
  \end{equation}
  Here 
  $$
  \cov{\zeta}{\tilde\zeta}
  :=\langle(\zeta-\langle\zeta\rangle)(\tilde\zeta-\langle\tilde\zeta\rangle)\rangle
  $$
  denotes the covariance of $\zeta$ and $\tilde\zeta$. Furthermore, for a site
  $y\in\mathbb{T}_L$, the random variable $\frac{\partial\zeta}{\partial y}$ is defined by
  $$
  (\frac{\partial\zeta}{\partial y})(\aa):=\zeta(\aa)-\int_{\Omega_0}\zeta(\aa)\,\beta(d\aa(y)).
  $$ 
\end{lemma}

\noindent
Like the classical partial derivatives $\{\frac{\partial\zeta}{\partial \aa_{kk}(y)}\}_{k=1,\cdots,d}$,
the function $\frac{\partial\zeta}{\partial y}$  measures how sensitively $\zeta$ 
depends on the variable $\aa(y)=\{\aa_{kk}(y)\}_{k=1,\cdots,d}\in\Omega_0$. 
For reasons explained in \cite{Gloria-Neukamm-Otto-11a} we call these derivatives {\it vertical}.
For obvious reasons, $\frac{\partial\zeta}{\partial y}$ is called the {\it discrete}
vertical derivative of $\zeta$ at $y$.

\medskip

\noindent
We note that for $\zeta=\tilde\zeta$, (\ref{O.1}) turns into the Spectral Gap Estimate (SG)
with constant 1, i.\ e.\
\begin{equation}\label{O.2}
  \langle\zeta^2\rangle\le\langle\sum_{y\in\mathbb{T}_L}(\frac{\partial\zeta}{\partial y})^2\rangle
\end{equation}
for any random variable $\zeta$ with vanishing expectation
$\langle\zeta\rangle=0$. We refer to \cite[Lemma~7]{Gloria-Neukamm-Otto-11a} for a proof 
in the present context and note that (\ref{O.2}) is extensively
used in the proof of Lemma~\ref{L:GNO1-1} as only probabilistic ingredient. 
This estimate is reminiscent of the Brascamp-Lieb inequality used by Naddaf and Spencer in \cite{Naddaf-Spencer-98}.
In a slightly different context, a covariance
estimate like (\ref{O.1}) was established in \cite[Lemma~3]{Gloria-Otto-09b}. 
For the convenience of the reader, we present the elementary proof of Lemma \ref{L1}.

\subsection{Corrector}\label{sec:corr}

We now introduce the important concept of the {\it corrector}. Since we only have to deal with
the {\it periodic} (as opposed to the infinite) ensemble $\langle\cdot\rangle$, 
we can avoid discussing all technicalities. Indeed, for any realization of $\aa$ according
to $\langle\cdot\rangle$, that is,
for any periodic coefficient field $\aa\in\Omega_L$ and for any coordinate direction
$j=1,\cdots,d$, there exists a unique scalar field $\phi_j(\aa,\cdot)\colon\T_L\rightarrow\mathbb{R}$
characterized by
\begin{equation}\label{eq:per-corr}
\nabla^*\aa(\nabla\phi_j+\ee_j)=0\;\;\mbox{on}\;\mathbb{T}_L\quad\mbox{and}\quad
\sum_{\mathbb{T}_L}\phi_j=0.
\end{equation}
Here $\ee_j$ denotes the unit vector in direction of the $j$-th coordinate axis, and the discrete derivatives are defined in \eqref{eq:def-disder}.
Clearly, for every $j=1,\cdots,d$, this defines a {\it random} scalar field $\phi_j$.
Evidently, this random field is {\it stationary} in the sense that 
for any shift $y\in\mathbb{R}^d$ one has $\phi_j(\aa(\cdot+y),\cdot)=\phi_j(\aa,\cdot+y)$.
The periodic function $\phi_j:\Z^d\rightarrow\R$ 
``corrects'' the affine function $x\mapsto x_j$ such that the resulting function
$\mathbb{Z}^d\ni x\mapsto \phi_j(x)+x_j$ is $a$-harmonic. In this sense,
$(\phi_1,\cdots,\phi_d)$ provide $\aa$-harmonic coordinates for $\mathbb{T}_L$.
We thus call the $\phi_j$'s the (periodic) {\it corrector}.
\medskip

\noindent
A crucial ingredient in the proof of our main result is the following
boundedness estimate on the moments of the corrector:
\begin{lemma}(Gloria, Neukamm \& Otto,
  \cite{Gloria-Neukamm-Otto-11a}).\label{L:GNO1-1}
  For $j=1,\ldots,d$ we have
  \begin{align}\label{eq:unif-bound-corr1}
    \expec{|\phi_{j}|^2}\,&\lesssim\, \left\{
      \begin{array}{ll}
        \ln L& \mbox{ for }d=2\\
        1&\mbox{ for }d>2
      \end{array}
    \right.,\\
    \label{eq:unif-bound-corr2}
    \expec{|\nabla\phi_{j}|^4}\,&\lesssim\, 1.
  \end{align}
  The multiplicative constants in these estimates only depend on
  $\lambda$ and $d$.
\end{lemma}
\noindent
We remark that the previous estimate indeed holds for the more general
class of stationary ensembles which satisfy SG. For the question in which sense $\lim_{L\uparrow\infty}\phi_j$
exists, which is however is not relevant for this paper, we refer to
\cite[Chapter~6]{Gloria-Neukamm-Otto-11a}. 

\medskip

\noindent
\subsection{Homogenized coefficient}
Let us recall that a coefficient field $$\aa\in\Omega:=\{\aa:\Z^d\rightarrow\Omega_0\}
=\Omega_0^{\mathbb{Z}^d}$$ can be seen as a description of a network of resistors: Suppose
$u\colon\mathbb{Z}^d\rightarrow\mathbb{R}$ is an $\aa$-harmonic function; if
$u(x)$ is interpreted as the potential at vertex $x$ and $\aa_{kk}(x)$ as conductivity
along the edge joining $x$ to $\ee_k+x$, then $\aa_{kk}(x)\nabla_ku(x)$ can be interpreted
as the (stationary) current along this edge. In this sense, the 
merit of the homogenized coefficient $\aa_{\hom}$ is that, almost surely, it relates the
spatially averaged potential gradient $\xi=\sum\nabla u$ to the spatially
averaged current $\aa_{\hom}\xi=\sum \aa\nabla u$. For our special $\aa$-harmonic function $u=\phi_j+x_j$,
the spatially averaged (over one period cell) potential gradient is given by
$$
L^{-d}\sum_{\mathbb{T}_L}(\nabla\phi_j+e_j)=e_j,
$$
whereas the spatially averaged current is given by
$$
L^{-d}\sum_{\mathbb{T}_L}\aa(\nabla\phi_j+e_j).
$$
Hence the expected value of the latter, that is,
\begin{equation}
\langle L^{-d}\sum_{\mathbb{T}_L}\aa(\nabla\phi_j+e_j)\rangle
\stackrel{\mbox{stationarity}}{=}\langle \aa(0)(\nabla\phi_j(0)+e_j)\rangle=:\aa_{\ho,L}e_j\label{eq:def-ahoml}
\end{equation}
is a good (and computable) proxy for $\aa_{\hom}$ as $L\uparrow\infty$. Qualitative homogenization
theory ensures that it is indeed true that
\begin{equation}
  \lim_{L\uparrow\infty}\aa_{\hom,L}=\aa_{\hom},\label{eq:def-ahom}
\end{equation}
see for instance \cite[Theorem~4.6]{Owhadi-03}. 
For us, this has the convenient side-effect
that we don't have to give the technically more demanding, independent definition
of $\aa_{\hom}$. For the latter, we refer to \cite[Introduction]{Gloria-Neukamm-Otto-11a} for instance. Moreover
--- and this is a second important ingredient for our result --- the quantitative
theory in \cite{Gloria-Neukamm-Otto-11a} provides an optimal estimate of this ``systematic error''
in the case of the (infinite) i.\ i.\ d.\ ensemble.

\begin{lemma}(Gloria, Neukamm \& Otto, \cite[Proposition~5]{Gloria-Neukamm-Otto-11a}).\label{L:GNO1-2}
Consider the i.\ i.\ d.\ ensemble. Then we have
$$
|\aa_{\hom,L}-\aa_{\hom}|\,\lesssim\, L^{-d}\ln^{d} L,
$$
where the multiplicative constant depends only on $d$ and $\lambda$.
\end{lemma}

\noindent
As opposed to the other important ingredients for our result, the annealed estimates
on the corrector and on the Green's function in Lemma~\ref{L:GNO1-1} and in Lemma~\ref{L:MO} below, 
the above ingredient relies on the i.\ i.\ d.\ property in a subtle way. 
In fact, if we were not dealing with the i.\ i.\ d.\ ensemble $\langle\cdot\rangle$, 
but with a more general infinite ensemble (that we call $\langle\cdot\rangle$ for the
purpose of the discussion in this paragraph only),
the choice of the periodic ensemble $\langle\cdot\rangle_L$
(that we endow with the index $L$ for the purpose of this discussion only) is more subtle
 --- and of practical importance: 
On the one hand, in view of Lemma \ref{L:GNO1-1} and of Lemma \ref{L:MO} below in particular,
$\langle\cdot\rangle_L$ should satisfy a Logarithmic Sobolev Inequality (LSI) uniformly in $L\uparrow\infty$.
On the other hand, in view of Lemma \ref{L:GNO1-2}, $\langle\cdot\rangle_L$ should be well-coupled to 
$\langle\cdot\rangle$ in order to avoid a large systematic error. If the ensemble is
not i.\ i.\ d.\ , the definition of $\langle\cdot\rangle_L$ through a ``brutal'' periodization, 
which is obtained by {\it restricting} $\aa(\cdot+y)$ to $\mathbb{T}_L$, with a random shift 
$y\in\mathbb{T}_L$ to retain stationarity, seems both unnatural and difficult to control.
It seems more natural and promising to us to define $\langle\cdot\rangle_L$ as the distribution
of $\aa$ under $\langle\cdot\rangle$ {\it conditioned} on the $\mathbb{T}_L$-periodicity of $\aa$.
However, this conditioning is singular; and we only expect control if 
the ensemble $\langle\cdot\rangle$ can be characterized by a sufficiently short-range 
(translation invariant) Hamiltonian. In this case, we expect 
that Dobrushin-Shlosman criteria (which as uniform mixing conditions ensure near-independence 
of $\aa(x)$ and $\aa(y)$ for $|x-y|\gg 1$ for all conditional measures) 
conveniently provide LSI for $\langle\cdot\rangle_L$ uniform in $L$.
The extension of Lemma~\ref{L:GNO1-2} to such a situation is
investigated in a forthcoming work, see \cite{GNO2}.
Let us also point out that the present proof relies on the covariance estimate of  Lemma~\ref{L1}. It is not yet clear to us whether such a covariance estimate can survive beyond the i.~i.~d. case.

\section{Proofs of Theorem~\ref{th:main} and Corollary~\ref{cor:1}}
\noindent
In this section we present the proofs of Theorem~\ref{th:main} and Corollary~\ref{cor:1}. Beforehand, we recall some auxiliary estimates on the elliptic Green's function that we need in the proof.

\subsection*{Structure of the proof and auxiliary estimates on the Green's function} The starting point to prove Theorem~\ref{th:main} is the same as in
the proof of \cite[Theorem~3]{Papanicolaou-Varadhan-79} by
Papanicolaou and Varadhan.
Recall that $z$ is given by \eqref{eq:uzhat}.
By uniform ellipticity and $L$-periodicity of $\aa$, an integration by parts on $\T_L$ yields
\begin{equation}\label{eq:GA1}
  \lambda \expec{\sum_{\T_L}|\nabla z|^2}\,\leq \, \expec{\sum_{\T_L}
    \nabla z\,\cdot\aa\nabla z}\,=\,\expec{\sum_{\T_L}
    z\, 
    \nabla^*\aa\nabla z}.
\end{equation}
As we shall see below,  in \textit{Step~1} of the proof of
Theorem~\ref{th:main}, an application of  $\nabla^*\aa\nabla$ to $z$ yields the decomposition
\begin{equation}\label{eq:rep-fo}
  \nabla^*\aa\nabla z\,=\,\nabla^* g+r_1+r_2,
\end{equation}
where the random vector field
$g:\T_L\to\R^d$, the deterministic scalar field $r_1:\T_L\to\R$, and the random
scalar field $r_2:\T_L\to\R$ are given by
\begin{align}\label{eq:12}  
\begin{split}
g_i\,=\,& -\sum_{j=1}^d\aa^{ii}\,\phi_j(\cdot+\ee_i)\,\nabla_i\nabla_ju_0,\\
r_1\,=\,& (\aa_{\ho,L}-\aa_{\ho}):\nabla^*\nabla u_0,\\
r_2\,=\,& (\bb-\aa_{\ho,L}):\nabla^*\nabla u_0.
\end{split}
\end{align}
Above, $\aa_{\ho,L}$ is defined via \eqref{eq:def-ahoml}, $\aa_{\hom}$ is defined via \eqref{eq:def-ahom}, and $\bb:\T_L\to\R^{d\times d}$ denotes the matrix field with entries
\begin{equation}\label{eq:def-bb}
  \bb^{ij}:=\aa^{ii}(\cdot-\ee_i)\big(\nabla_i\phi_{j}(\cdot-\ee_i)+\delta(i-j)\big)\qquad\text{for }i,j=1,\ldots,d.
\end{equation}
As we shall see below (in \textit{Step~3} of the proof of Theorem~\ref{th:main}), the scalar fields $r_1$ and $r_2$ satisfy
\begin{equation}\label{eq:14}
  \sum_{\T_L}r_1=0\qquad\mbox{and}\qquad\expec{r_2(x)}=0\text{  for all $x\in\T_L$.}
\end{equation}
By combining \eqref{eq:GA1} and \eqref{eq:rep-fo}, the bound \eqref{eq:mainL} follows from estimates of the terms $\expec{\sum_{\T_L} \nabla z\cdot g}$, $\expec{\sum_{\T_L} (z-\bar z)\,r_1}$ and $\expec{\sum_{\T_L} z\,r_2}$, an estimate of the spatial mean
\begin{equation*}
  \bar z:=L^{-d}\sum_{\T_L}z,
\end{equation*}
and the discrete Poincar\'e inequality on the torus.

\smallskip

\noindent
The most intricate estimate is the one of $\expec{\sum_{\T_L}
  z\,r_2}$, see \textit{Step~5} in the proof of
Theorem~\ref{th:main} below.
Since the term $r_2(x)$ has vanishing expectation for all $x\in\T_L$ one may write it as a covariance
$$
\expec{\sum_{x\in \T_L} z(x)r_2(x)}\,=\,\sum_{x\in \T_L} \cov{z(x)}{r_2(x)}.
$$
In order to benefit from this, we appeal to the covariance estimate of Lemma~\ref{L1} and to a vertical derivative calculus on coefficient fields that we introduced in \cite{Gloria-Neukamm-Otto-11a}.
The process of estimating the vertical derivatives of $r_2$ and $z$ involves the (periodic) Green's function:
\begin{defi}\label{lem:def-GreenL}
  The $L$-periodic Green's function $G_L:\T_L\times \T_L\times
  \Omega_L \to \R$ is defined as follows. For all $y\in \T_L$ and $\aa\in\Omega_L$
  the function $G_L(\cdot,y;\aa)$ is the unique $L$-periodic mean free solution to
  \begin{equation}\label{eq:def-GreenL}
    \nabla^*\aa\nabla G_L(\cdot,y;\aa)\,=\,\delta(\cdot-y)-L^{-d} \qquad \text{on }\T_L,
  \end{equation}
where $\delta$ is the Dirac mass at zero.
\end{defi}
\noindent When no confusion occurs we use the shorthand notation $G_L(x,y)$ for $G_L(x,y;\aa)$.
We shall use both quenched (i.~e. pointwise deterministic)  and annealed (i.~e. statistically averaged) estimates on $|\nabla G_L|$.
The pointwise estimates rely on the \mbox{De
  Giorgi}-Nash-Moser H\"older regularity theory (and are standard in the continuum case):
\begin{lemma}\label{lem:quenched}
There exists $\gamma>0$ depending only on $\lambda$ and $d$ such that for all $\aa\in \Omega_L$ and $L\in \N$, the Green's function $G_L(\cdot,\cdot;\aa)$ satisfies the following quenched estimate:
\begin{equation}\label{eq:quenched}
|\nabla_x G_L(x,y;\aa)|,|\nabla_y G_L(x,y;\aa)|\,\lesssim\, (|x-y\mod L|+1)^{2-d-\gamma}.
\end{equation}
In the estimate the multiplicative constant only depends on $\lambda$ and $d$.
(Note that $G_L$ is symmetric, so that the estimate \eqref{eq:quenched} does not depend on the variable with respect to which we differentiate.)
\end{lemma}
\noindent 
See Appendix~\ref{sec:proof2} for the proof. The crucial other ingredient is the recent annealed estimate of \cite{Marahrens-Otto-13} by
Marahrens and the third author, which we recall below in a version for
the $L$-periodic Green's functions:
\begin{lemma}\label{L:MO}(Marahrens \& Otto,\cite[Theorem~1]{Marahrens-Otto-13}).
The periodic Green's function $G_L$ satisfies the following annealed estimates:
\begin{equation}\label{eq:annealed}
  \expec{|\nabla G_L(x,y)|^{4}}^{\frac{1}{4}}\,\lesssim\, (|x-y\mod L|+1)^{1-d},
\end{equation}
and 
\begin{equation}\label{eq:annealed2}
\expec{|\nabla_x\nabla_{y}  G_L(x,y)|^{4}}^{\frac{1}{4}}\,\lesssim\, (|x-y\mod L|+1)^{-d},
\end{equation}
where the multiplicative constants only depend on $\lambda$ and $d$.
\end{lemma}
\medskip
\noindent Let us mention that the proof of the annealed estimates on the derivatives of the Green's functions relies on a strengthened version of the spectral gap estimate \eqref{O.2}, namely a Logarithmic-Sobolev Inequality.
We refer the reader to \cite{Marahrens-Otto-13} for details.
Note that as for the variance estimate by Naddaf and Spencer in \cite{Naddaf-Spencer-98}, an optimal control of the fourth moment is enough for our quantitative expansion.


\subsection*{Proof of Theorem~\ref{th:main}}
\noindent
We prove the estimate of Theorem~\ref{th:main} in its rescaled formulation
\eqref{eq:mainL}.  Note that the identity
\begin{equation}\label{eq:discrete-hessian}
  \sum_{\T_L}\sum_{i,j=1}^d|\nabla^*_i\nabla_j u_0|^2=\sum_{\T_L}\sum_{i,j=1}^d|\nabla_i\nabla_j u_0|^2=\sum_{\T_L}|\nabla^2 u_0|^2
\end{equation}
follows from periodicity and the elementary identity $|\nabla^*_i\nabla_ju_0(\cdot)|=|\nabla_i\nabla_ju_0(\cdot-\ee_i)|$.
The argument is divided into five steps. In the
first step we derive the decomposition  \eqref{eq:rep-fo} with property \eqref{eq:14}.
In Steps~2 and~3 we argue that it suffices to prove
\begin{equation}\label{pf:th:main}
  \expec{\sum_{\T_L}|\nabla z|^2}\lesssim\,\mu_d(L)\,\sum_{\T_L}|\nabla^2 u_0|^2.
\end{equation}
In the remaining steps we prove \eqref{pf:th:main} starting with the inequality
$\lambda\expec{\sum_{\T_L} |\nabla z|^2} \, \leq \,
\expec{\sum_{\T_L}z\,\nabla^*\aa\nabla z}$. To that end we appeal to the representation formula~\eqref{eq:rep-fo} which is a sum of a term in
divergence form $\nabla^*g$, a deterministic term $r_1$ and a remainder
with vanishing expectation $r_2$. The first two terms are estimated in
Step~4. The third term $r_2$, which has vanishing expectation, is
controlled using the covariance estimate of
Lemma~\ref{L1} and will be treated in Steps 5a--5d. 
\medskip

\step 1 Derivation of the decomposition \eqref{eq:rep-fo} with property \eqref{eq:14}.

\medskip

\noindent
Let us show that \eqref{eq:rep-fo} holds with $g,r_1$ and $r_2$ given by \eqref{eq:12}. By the defining equation~\eqref{eq:torus-L} for $u $,
\begin{eqnarray}
  \nabla^*\aa\nabla z&\stackrel{\eqref{eq:uzhat}}{=}&\nabla^*\aa\nabla u -\nabla^*\aa\nabla
  u_0-\sum_{j=1}^d\nabla^*\aa\nabla(\phi_{j}\nabla_j u_0) \nonumber \\
  &\stackrel{\eqref{eq:torus-L}}{=}& \tilde f -\underbrace{\nabla^*\aa\nabla
  u_0}_{\dps =:I}-\underbrace{\sum_{j=1}^d\nabla^*\aa\nabla(\phi_{j}\nabla_j u_0)}_{\dps =:II}.\label{eq:prth-2}
\end{eqnarray}
We shall use the following discrete Leibniz rule: For all $\zeta_1,\zeta_2:\T_L\to \R$,
\begin{align}\label{eq:leibniz}
    \nabla_{i}(\zeta_1\zeta_2)\,=\,(\nabla_i\zeta_1)\,\zeta_2+\zeta_1(\cdot+\ee_i)\nabla_i\zeta_2,\quad
    \nabla_{i}^*(\zeta_1\zeta_2)\,=\,(\nabla^*_{i}
    \zeta_1)\,\zeta_2+\zeta_1(\cdot-\ee_i)\nabla^*_i\zeta_2.
\end{align}
For the first term this yields
\begin{equation*}
  I\,=\,\sum_{i=1}^d(\nabla_i^*\aa^{ii})\,\nabla_{i}u_0+\sum_{i=1}^d \aa^{ii}(\cdot-\ee_i)\,\nabla^*_{i}\nabla_iu_0,
\end{equation*}
while for the second term we obtain
\begin{eqnarray*}
II&=&\sum_{i,j=1}^d\nabla^*_i\Big(\, \aa^{ii}\nabla_i(\phi_j\nabla_ju_0)\,\Big)\\
&\stackrel{\eqref{eq:leibniz}}{=}&\sum_{i,j=1}^d\nabla^*_i\Big(\,
\aa^{ii}(\nabla_i\phi_j\nabla_ju_0+\phi_j(\cdot+\ee_i)\nabla_i\nabla_ju_0)\,\Big)\\
&\stackrel{\eqref{eq:leibniz}}{=}&\sum_{i,j=1}^d\nabla^*_i(\,
\aa^{ii}\nabla_i\phi_j\,)\nabla_ju_0\ +\ \sum_{i,j=1}^d\aa^{ii}(\cdot-\ee_i)\nabla_i\phi_j(\cdot-\ee_i)\nabla^*_i\nabla_j u_0\\
&&+\sum_{i,j=1}^d\nabla^*_i\Big(\,
\aa^{ii}\,\phi_j(\cdot+\ee_i)\nabla_i\nabla_ju_0\,\Big).
\end{eqnarray*}
Because of the periodic corrector equation \eqref{eq:per-corr}, the first term of the r.~h.~s. turns into
\begin{equation*}
  \sum_{i,j=1}^d\nabla^*_i(\,
  \aa^{ii}\nabla_i\phi_j\,)\,\nabla_ju_0=-\sum_{j=1}^d\big(\,\nabla^*_j\aa^{jj}\,\big)\nabla_ju_0.
\end{equation*}
Hence, the terms in $I+II$ that involve the \emph{first} derivative
of $u_0$ cancel, so that
\begin{eqnarray}\notag
I+II&=&
\sum_{i,j=1}^d\aa^{ii}(\cdot-\ee_i)\big(\,\delta(i-j)+\nabla_i\phi_j(\cdot-\ee_i)\,\big)\nabla^*_i\nabla_j
u_0\\
&&+\sum_{i,j=1}^d\nabla^*_i\Big(\,
\aa^{ii}\,\phi_j(\cdot+\ee_i)\nabla_i\nabla_ju_0\,\Big)
.\label{eq:prth-1}
\end{eqnarray}
The last term on the r.~h.~s. is precisely $-\nabla^*g$, the
first term is $\bb:\nabla^2 u_0$, cf. \eqref{eq:12} and \eqref{eq:def-bb}.
The claim \eqref{eq:rep-fo} then follows from \eqref{eq:prth-2},
\eqref{eq:prth-1} and identity \eqref{eq:9b} which can be written in
the form $\tilde f\,=\,\aa_\ho:\nabla^*\nabla u_0$ since $\aa_\ho$ is constant.

\noindent
To conclude this step, we prove \eqref{eq:14}. The first identity
simply follows from the $L$-periodicity of $u_0$. The second
identity can be seen as follows: By the definition of $\bb$, the stationarity of $\aa$ and $\phi_{j}$, and the definition of $\aa_{\ho,L}$ via \eqref{eq:def-ahoml} we have
\begin{equation*}
\expec{\bb^{ij}(\cdot)}=
\expec{\aa^{ii}(0)\Big(\delta(i-j)+\nabla_{i}\phi_{j}(0)\Big)}\notag=\aa_{\hom,L}^{ij},
\end{equation*}
so that
\begin{eqnarray*}
  \expec{\bb:\nabla^2 u_0}=\expec{\bb}:\nabla^2 u_0=\aa_{\ho,L}:\nabla^2 u_0
\end{eqnarray*}
as desired.

\medskip

\step 2 Reduction to an estimate for $\nabla z$.

\medskip

\noindent
We claim that \eqref{eq:mainL} (and thus the statement of
Theorem~\ref{th:main}, see Remark~\ref{R:rescaling})
follows from \eqref{pf:th:main}. Indeed, by the discrete Poincar\'e
inequality $\sum_{\T_L}z^2\lesssim L^2\sum_{\T_L}|\nabla
z|^2+L^d\bar z^2$ we only need to prove that
\begin{equation*}
  \expec{\bar z^2}\lesssim L^{2-d}\mu_d(L)\,\sum_{\T_L}|\nabla^2u_0|^2.
\end{equation*}
Since the spatial means of $u $ and $u_0$ vanish by definition, we have $\bar
z=-L^{-d}\sum_{\T_L}\sum_{j=1}^d\phi_j\nabla_ju_0$,
so that
\begin{eqnarray*}
  \expec{\bar
    z^2}&=&L^{-2d}\sum_{x\in\T_L}\sum_{x'\in\T_L}\expec{\big(\,\sum_{i=1}^d\phi_i(x)\nabla_iu_0(x)\,\big)\big(\,\sum_{j=1}^d\phi_j(x')\nabla_ju_0(x')\,\big)}.
\end{eqnarray*}
We expand the square on the r.~h.~s.. Since $u_0$ is
deterministic we get
\begin{equation*}
  \expec{\bar
    z^2}
  \,=\, L^{-2d}\sum_{x\in\T_L}\sum_{x'\in\T_L}\sum_{i,j=1}^d\nabla_iu_0(x)\nabla_ju_0(x')\expec{\phi_i(x)\phi_j(x')}.
\end{equation*}
By Cauchy-Schwarz' inequality, stationarity of the correctors, and the bounds of Lemma~\ref{L:GNO1-1},
$$
\expec{\phi_i(x)\phi_j(x')}\,\leq\, \max_{k=1,\ldots,d} \expec{\phi_k^2} \,\lesssim \, \mu_d(L).
$$
Hence, 
\begin{equation*}
  \expec{\bar z^2}  \,\lesssim \, \mu_d(L)  \left(L^{-d}\sum_{x\in\T_L}|\nabla
  u_0|\right)^2.
\end{equation*}
The desired estimate then follows from Jensen's and Poincar\'e's inequalities.

\medskip

\step 3 Reduction based on the decomposition \eqref{eq:rep-fo}.

\medskip

\noindent
In this step we argue that the desired estimate reduces to the following statement: Suppose that
the functions $g$, $r_1$ and $r_2$ of the decompositon \eqref{eq:rep-fo} satisfy the following estimates in addition to \eqref{eq:14}:
\begin{eqnarray}
  \label{pf:main:est1}
  \expec{\sum_{\T_L}|g|^2}&\lesssim& \mu_d(L)\,\sum_{\T_L}|\nabla^2 u_0|^2,\\
  \label{pf:main:est2}
 \sum_{\T_L}r_1^2&\lesssim&
  \tfrac{\ln^d L}{L^d}\,\sum_{\T_L}|\nabla^2 u_0|^2,\\
  \label{pf:main:est3}
  \expec{\sum_{\T_L}z\,r_2}&\lesssim&\expec{\sum_{\T_L}|\nabla
    z|^2}^{\frac{1}{2}}\left(\sum_{\T_L}|\nabla^2 u_0|^2\right)^{\frac{1}{2}}\,+\,\mu_d(L)\sum_{\T_L}|\nabla^2 u_0|^2,
\end{eqnarray}
then estimate \eqref{eq:mainL} (and thus Theorem~\ref{th:main}) holds.

\smallskip

\noindent
By Step~2 we just have to check \eqref{pf:th:main}.
Indeed, by combining \eqref{eq:GA1}, the decomposition \eqref{eq:rep-fo}, the triangle inequality, an integration
by parts, and \eqref{eq:14} we get
\begin{eqnarray*}
  \lambda\expec{\sum_{\T_L}|\nabla z|^2}&\leq&
  \left|\expec{\sum_{\T_L}\nabla z\cdot g}\right| +  \left|\expec{\sum_{\T_L}(z-\bar z)r_1}\right| +  \left|\expec{\sum_{\T_L}z\,r_2}\right|.
\end{eqnarray*}
The first term is estimated by \eqref{pf:main:est1} and the
Cauchy-Schwarz inequality:
\begin{eqnarray*}
  \left|\expec{\sum_{\T_L}\nabla z\cdot g}\right| &\lesssim& \mu_d^{\frac{1}{2}}(L)\,\expec{\sum_{\T_L}|\nabla
    z|^2}^{\frac{1}{2}}\,\left(\sum_{\T_L}|\nabla^2 u_0|^2\right)^{\frac{1}{2}}.
\end{eqnarray*}
The second term is estimated by \eqref{pf:main:est2}, the
Cauchy-Schwarz inequality, the Poincar\'e inequality for functions
on $\T_L$ with zero mean, and the elementary estimate  $\frac{\ln^dL}{L^{d-1}}\lesssim 1$ for $d>1$:
\begin{eqnarray*}
  \left|\expec{\sum_{\T_L}(z-\bar z) r_1}\right| &\lesssim& \expec{\sum_{\T_L}|\nabla
    z|^2}^{\frac{1}{2}}\,\left(\sum_{\T_L}|\nabla^2 u_0|^2\right)^{\frac{1}{2}}.
\end{eqnarray*}
The combination of \eqref{pf:main:est3} with the previous three
inequalities yields
\begin{equation*}
  \expec{\sum_{\T_L}|\nabla z|^2}\,\lesssim\,   \mu_d^{\frac{1}{2}}(L)\,\expec{\sum_{\T_L}|\nabla
    z|^2}^{\frac{1}{2}}\,\left(\sum_{\T_L}|\nabla^2 u_0|^2\right)^{\frac{1}{2}}+ \mu_d(L)\,\sum_{\T_L}|\nabla^2 u_0|^2  ,
\end{equation*}
which implies \eqref{pf:th:main} by Young's inequality.

\medskip

\step 4 Proof of the estimates \eqref{pf:main:est1} and
\eqref{pf:main:est2}.
\medskip

\noindent
Estimate \eqref{pf:main:est1} follows from the definition of $g$, cf. \eqref{eq:12}, the
bound \eqref{eq:unif-bound-corr1} of Lemma~\ref{L:GNO1-1}
on the second moment of the stationary $\phi_{j}$ and identity \eqref{eq:discrete-hessian}. Similarly, \eqref{pf:main:est2} follows from the definition of $r_1$ and the optimal bound on the error $|\aa_{\ho,L}-\aa_{\ho}|$ of
Lemma~\ref{L:GNO1-2}.

\medskip

\noindent In the last step we prove \eqref{pf:main:est3}.
Since this step is rather long, we subdivide it further.

\medskip

\step{5a} Application of the
covariance estimate.

\medskip

\noindent 
Since by \eqref{eq:14}, $\expec{r_2(x)}=0$ for all $x\in \T_L$, we have
\begin{align*}
  &\left|\expec{\sum_{x\in\T_L}z(x)r_2(x)}\right|\,=\,\left|\sum_{x\in\T_L}
  \cov{z(x)}{r_2(x)}\right|\\
  &\qquad\stackrel{\eqref{eq:12}}{=}\,
  \left|\sum_{x\in\T_L}\sum_{i,j=1}^d\nabla_i^*\nabla_j
  u_0(x)\,\cov{z(x)}{\aa_{\ho,L}^{ij}-\bb^{ij}(x)}\right|\\
  &\qquad\stackrel{\eqref{eq:discrete-hessian}}{\leq}\,
  \left(\sum_{\T_L}|\nabla^2 u_0|^2\right)^{\frac{1}{2}}\,
  \left(\sum_{x\in\T_L}\sum_{i,j=1}^d\left(\cov{z(x)}{\aa_{\ho,L}^{ij}-\bb^{ij}(x)}\right)^2\right)^{\frac{1}{2}}.
\end{align*}
With the covariance
estimate in Lemma~\ref{L1} the r.~h.~s. is bounded by
\begin{equation*}
  \left(\sum_{\T_L}|\nabla^2u_0|^2\right)^{\frac{1}{2}}\,\left(\sum_{x\in\T_L}
  \sum_{i,j=1}^d\left(\sum_{y\in\T_L}\expec{\left(\frac{\partial
          z(x)}{\partial y}
\right)^2}^{\frac{1}{2}}\expec{\left(\,\frac{\partial \bb^{ij}(x)}{\partial y}\,\right)^2}^{\frac{1}{2}}\right)^2\right)^{\frac{1}{2}}.
\end{equation*}
Hence, for \eqref{pf:main:est3} it suffices to prove that
\begin{align}\label{pf:main:est3-alt}
  \begin{split}
    &\left(\sum_{x\in\T_L}\left(\sum_{y\in\T_L}\expec{\left(\frac{\partial
            z(x)}{\partial y}
        \right)^2}^{\frac{1}{2}}\expec{\left(\,\frac{\partial
            \bb^{ij}(x)}{\partial
            y}\,\right)^2}^{\frac{1}{2}}\right)^2\right)^{\frac{1}{2}}\\
    &\qquad\leq \expec{\sum_{\T_L}|\nabla
      z|^2}^{\frac{1}{2}}\,+\,\mu_d(L)\,\left(\sum_{\T_L}|\nabla^2 u_0|^2\right)^{\frac{1}{2}}\qquad\text{for }i,j=1,\ldots,d.
  \end{split}
\end{align}
To estimate the vertical derivatives we need to identify $\frac{\partial z}{\partial y}$ and $\frac{\partial \bb^{ij}}{\partial
  y}$. This is done by appealing to the elliptic equations \eqref{eq:9} and
\eqref{eq:per-corr} and vertical differential calculus. Since the
basic argument is simple, but polluted due to the discrete nature of the vertical
and spatial derivatives, we will first present a \emph{formal} calculation where the vertical derivative $\frac{\partial}{\partial
  y}$ is replaced by the classical partial derivative
$\frac{\partial}{\partial\aa_{kk}(y)}$ (defined for differentiable functions on $\Omega_L$). 
The \emph{rigorous} argument is then carried out in Step~5c and Step~5d below.

\medskip

\step{5b} \emph{Formal} derivation of formulas for the vertical derivatives.

\medskip

\noindent 
We first (formally) identify $\frac{\partial u }{\partial\aa_{kk}(y)}$ and
$\frac{\partial\phi_i(x)}{\partial\aa(x)}$. Applying $\frac{\partial}{\partial\aa_{kk}(y)}$ to the elliptic equations
\eqref{eq:9} and \eqref{eq:per-corr} yields
\begin{eqnarray*}
  \nabla^*\aa(x)\nabla\frac{\partial
    u(x)}{\partial\aa_{kk}(y)}&=&-\nabla^*\frac{\partial\aa(x)}{\partial\aa_{kk}(y)}\nabla
  u(x),\\
  \nabla^*\aa(x)\nabla\frac{\partial
    \phi_j(x)}{\partial\aa_{kk}(y)}&=&-\nabla^*\frac{\partial\aa(x)}{\partial\aa_{kk}(y)}(\nabla
  \phi_j(x)+\ee_j),
\end{eqnarray*}
using that $\frac{\partial}{\partial\aa_{kk}(y)}$ and $\nabla$
commute. Since
$\frac{\partial\aa(x)}{\partial\aa_{kk}(y)}=(\ee_k\otimes\ee_k)\delta(x-y\mod
L)$ for all $x,y\in\T_L$,
the Green representation formula yields
\begin{align}\label{eq:1}
  \begin{split}
    \frac{\partial
      u(x)}{\partial\aa_{kk}(y)}\,=&\,-\nabla_{y_k}G_L(x,y)\nabla_{k}
    u(y),\\
    \frac{\partial
      \phi_j(x)}{\partial\aa_{kk}(y)}\,=&\,-\nabla_{y_k}G_L(x,y)(\nabla_k
    \phi_j(y)+\delta(k-j)).
  \end{split}
\end{align}
Next we identify $\frac{\partial\bb^{ij}}{\partial\aa_{kk}(y)}$. We
apply $\frac{\partial}{\partial\aa_{kk}(y)}$ to the definition
\eqref{eq:def-bb} of $\bb^{ij}$ and use the identity above in the form of
\begin{equation}\label{eq:5}
  \frac{\partial\nabla_i\phi_j(x)}{\partial\aa_{kk}(y)}\,=\,-\nabla_{x_i}\nabla_{y_k}G_L(x,y)(\nabla_k\phi_j(y)+\delta(k-j)).
\end{equation}
Rearranging the terms yields the identity
\begin{align}\label{eq:pd-bb}
  \begin{split}
    \frac{\partial\bb^{ij}(x+\ee_i)}{\partial\aa_{kk}(y)}\,=\,&\delta(k-i)\delta(x-y\mod
    L)(\nabla_i\phi_j(y)+\delta(i-j))\\
    &-\aa^{ii}(x)\nabla_{x_i}\nabla_{y_k}G_L(x,y)(\nabla_k\phi_j(y)+\delta(k-j)).
  \end{split}
\end{align}
Likewise, for the identification of  $\frac{\partial
  z(x)}{\partial\aa_{kk}(y)}$ we apply $\frac{\partial}{\partial\aa_{kk}(y)}$ to \eqref{eq:uzhat}:
\begin{eqnarray*}
  \frac{\partial z(x)}{\partial\aa_{kk}(y)}&=&
  \frac{\partial u(x)}{\partial\aa_{kk}(y)}
  -\sum_{j=1}^d\frac{\partial\phi_j(x)}{\partial\aa_{kk}(y)}\nabla_ju_0(x)\\
  &\stackrel{\eqref{eq:1}}{=}&
  -\nabla_{y_k}G_L(x,y)\left(\nabla_k u(y) - \sum_{j=1}^d(\nabla_k\phi_j(y)+\delta(k-j))\nabla_ju_0(x)\right).
\end{eqnarray*}
Since we want to make $\nabla z$ appear, we
substitute $\nabla u$ by the following expression
\begin{equation}\label{eq:6}
  \nabla u(y)=\nabla z(y)+\sum_{j=1}^d(\nabla\phi_j(y)+\ee_j)\nabla_ju_0(y)+\sum_{j=1}^d\phi_j(y)\nabla\nabla_j u_0(y),
\end{equation}
which can formally be obtained by applying $\nabla$ to
\eqref{eq:uzhat} and using the continuum Leibniz rule
$\nabla(\phi_j\nabla_j u_0)=\nabla\phi_j\nabla_j
u_0+\phi_j\nabla\nabla_j u$. We then get
\begin{align}\label{eq:pd-z}
  \begin{split}
    &\frac{\partial z(x)}{\partial\aa_{kk}(y)}=
    -\nabla_{y_k}G_L(x,y)\Bigg(\,\nabla_k
    z(y)+\sum_{j=1}^d(\nabla_k\phi_j(y)+\delta(k-j))(\nabla_ju_0(y)-\nabla_ju_0(x))\\
    &\qquad\qquad\qquad+\sum_{j=1}^d\phi_j(y)\nabla_k\nabla_j
    u_0(y)\,\Bigg).
  \end{split}
\end{align}
Let us stress the fact that the expression in the brackets on the r.~h.~s. is
$\nabla z(y)$ plus terms that vanish if $u_0$ is affine. This will be
crucial in order to obtain an optimal estimate. 

\medskip

\step{5c} \emph{Rigorous} derivation of formulas for the vertical derivatives.

\medskip

\noindent  We now derive rigorous versions of
\eqref{eq:pd-bb} and \eqref{eq:pd-z}, which will lead to the desired
estimate \eqref{pf:main:est3-alt}. For the rigorous argument
$\frac{\partial}{\partial\aa_{kk}(y)}$ has to be replaced by the
discrete vertical derivative $\frac{\partial}{\partial y}$ for which the Leibniz rule is
not valid. The main ingredient is the following rigorous version of \eqref{eq:1}: For $j=1,\ldots,d$ and $x,y\in \T_L$ we have
\begin{align}\label{eq:vd-phi}
  \frac{\partial\phi_j(x)}{\partial y}\,&=\,-\nabla_yG_L(x,y)\cdot[\nabla\phi_j(y)+\ee_j]_y,\\
  \frac{\partial\nabla_i\phi_j(x)}{\partial y}\,&=\,-\nabla_{x_i}\nabla_yG_L(x,y)\cdot[\nabla\phi_j(y)+\ee_j]_y \label{eq:vd-grad-phi}\\
  \label{eq:vd-uL}
  \frac{\partial u(x)}{\partial
    y}\,&=\,-\nabla_yG_L(x,y)\cdot [\nabla u(y)]_y,
\end{align}
where $[\cdot]_y$ denotes the commutator of the multiplication with $\aa(y)$ and $\frac{\partial}{\partial y}$, i.~e.
\begin{equation}\label{eq:def-com}
  [F]_y:=\frac{\partial(\aa(y)F)}{\partial y}-\aa(y)\frac{\partial F}{\partial y}=\aa(y)\expec{F}_y-\expec{\aa(y)F}_y
\end{equation}
for all random vectors $F$.
By Jensen's inequality,  the commutator satisfies the following estimate: For all $1\leq q < \infty$,
\begin{equation}\label{referee}
\expec{ |[F]_y|^q} \, \leq \, 2^q \expec{|F|^q},
\end{equation}
which we will use in the sequel for $q=2$ and $q=4$.
Note that for all $x,y\in \T_L$, since the coefficients are i.~i.~d., 
\begin{equation}\label{eq:iid-commut}
  \frac{\partial(\aa(x)F)}{\partial y}-\aa(x)\frac{\partial F}{\partial y}=[F]_y\delta(x-y \mod L).
\end{equation}
Here comes the argument for \eqref{eq:vd-phi}. By \eqref{eq:per-corr} we have
\begin{eqnarray*}
  0&=&\nabla^*\aa(x)(\nabla\phi_j(x)+\ee_j)-\expec{\nabla^*\aa(x)(\nabla\phi_j(x)+\ee_j)}_{y}\\
  &=&\nabla^*\Big(\,\aa(x)(\nabla\phi_j(x)+\ee_j)-\expec{\aa(x)(\nabla\phi_j(x)+\ee_j)}_{y}\,\Big)\\
  &=&\nabla^*\Big(\,\aa(x)\frac{\partial \nabla\phi_j(x)}{\partial y}+\aa(x)\expec{\nabla\phi_j(x)+\ee_j}_{y}-\expec{\aa(x)(\nabla\phi_j(x)+\ee_j)}_{y}\,\Big) 
\end{eqnarray*}
Using then \eqref{eq:def-com} and \eqref{eq:iid-commut} for $F=\nabla\phi_j(x)+\ee_j$, this turns into
\begin{equation*}
\nabla^*\Big(\aa(x)\nabla \frac{\partial \phi_j(x)}{\partial y}\Big) \,=\,-\nabla^*\Big(   [\nabla\phi_j(x)+\ee_j]_y\delta(x-y \mod L)\Big),
\end{equation*}
from which \eqref{eq:vd-phi} follows by the Green representation formula and an integration by parts.
Identity \eqref{eq:vd-grad-phi} follows from applying $\nabla_x$ to \eqref{eq:vd-phi}. The argument
for \eqref{eq:vd-uL} is similar to the one for \eqref{eq:vd-phi} and
left to the reader.

\smallskip

\noindent We are now in position to derive the rigorous versions of \eqref{eq:pd-bb} and \eqref{eq:pd-z},
and start with $b^{ij}$. We claim that 
\begin{align}\label{eq:main-step4-1}
  \begin{split}
    \frac{\partial\bb^{ij}(x+\ee_i)}{\partial y}\,=\,\Big(\delta(x-y\mod L)\ee_i-\aa^{ii}(x)\nabla_{x_i}\nabla_y
    G_L(x,y)\Big)\cdot [\nabla\phi_j(y)+\ee_j]_y,
  \end{split}
\end{align}
Indeed, using again \eqref{eq:def-com} for $F=\nabla\phi_j(x)+\ee_j$, we have
\begin{eqnarray*}
  \frac{\partial(\aa(x)(\nabla\phi_j(x)+\ee_j))}{\partial
    y}&=&\aa(x)\nabla\frac{\partial\phi_j(x)}{\partial
    y}+\aa(x)\expec{\nabla\phi_j(x)+\ee_j}_{y}-\expec{\aa(x)(\nabla\phi_j(x)+\ee_j)}_{y}\\
  &=&\aa(x)\frac{\partial\nabla\phi_j(x)}{\partial
    y}+\delta(x-y\mod L)[\nabla\phi_j(y)+\ee_j]_y,
\end{eqnarray*}
which, combined with \eqref{eq:vd-grad-phi}, yields \eqref{eq:main-step4-1}.

\smallskip

\noindent We then turn to $z$ and claim that
  \begin{align}
    \begin{split}
      \label{eq:main-step4-2}
      \frac{\partial z(x)}{\partial y}\,=\,&-\nabla_yG_L(x,y)\cdot \Big( F_1(y) +F_2(y,x)+F_3(y)\Big) ,
    \end{split}
  \end{align}
  where
\begin{eqnarray*}
    F_1(y)&:=&[\nabla z(y)]_{y},\qquad F_2(y,x)\,:=\,\sum_{j=1}^d\big(\nabla_ju_0(y)-\nabla_ju_0(x)\big)[\nabla\phi_j(y)+\ee_j]_y,\\
    F_3(y)&:=&\sum_{i,j=1}^d\nabla_i\nabla_ju_0(y)[\phi_j(y+\ee_i)\ee_i]_y.
  \end{eqnarray*}
 We first apply
  $\frac{\partial}{\partial y}$ to \eqref{eq:uzhat} and use
  \eqref{eq:vd-phi} and \eqref{eq:vd-uL}:
  \begin{equation}\label{eq:main-step4-8}
    \frac{\partial z(x)}{\partial y}\,=\,-\nabla_yG_L(x,y)\cdot\Big([\nabla u(y)]_y-\sum_{j=1}^d\nabla_j
    u_0(x)\,[\nabla\phi_j(y)+\ee_j]_y\,\Big).
  \end{equation}
  We then wish to substitute the term $\nabla u(y)$ by
  an expression that involves $\nabla z(y)$. To that end we apply
  $\nabla$ to \eqref{eq:uzhat} and get with the help of the discrete
  Leibniz rule \eqref{eq:leibniz}
  \begin{eqnarray*}
    \nabla u(y)&=&\nabla z(y)+\nabla u_0(y)
    +\sum_{j=1}^d\nabla_ju_0(y)\nabla\phi_{j}(y) +\sum_{i,j=1}^d\nabla_i\nabla_ju_0(y)\,\phi_{j}(y+\ee_i)\ee_i\\
    &=&\nabla z(y)+\sum_{j=1}^d\nabla_ju_0(y)(\nabla\phi_j(y)+\ee_j)+\sum_{i,j=1}^d\nabla_i\nabla_ju_0(y)\,\phi_{j}(y+\ee_i)\ee_i.
  \end{eqnarray*}
Combined with \eqref{eq:main-step4-8} the desired identity \eqref{eq:main-step4-2} follows.

\medskip

\step{5d} Estimates of the vertical derivatives of $\bb^{ij}$ and $z$.

\medskip

\noindent
We claim that
\begin{equation}\label{pf:main-step5b}
    \expec{\left(\frac{\partial\bb^{ij}(x)}{\partial
          y}\right)^2}^{{1}/{2}}\,\lesssim\,(|y-x\mod L|+1)^{-d}
\end{equation}
and
 \begin{equation}
    \label{eq:7}
    \expec{\left(\frac{\partial z(x)}{\partial
          y}\right)^2}^{{1}/{2}}\,\qquad\lesssim I_1+I_2,
  \end{equation}
  where
  \begin{eqnarray*}
    I_1&:=&(|y-x \mod L|+1)^{2-d-\gamma}\left(\expec{|\nabla
      z(y)|^2}^{1/2}+\mu_d(L)|\nabla^2 u_0(y)|\right),\\
    I_2&:=&(|y-x \mod L|+1)^{1-d}|\nabla u_0(x)-\nabla u_0(y)|,
  \end{eqnarray*}
  for some $\gamma>0$ depending only on $\lambda$ and $d$. 

\smallskip

\noindent We start with \eqref{pf:main-step5b}.
By the Cauchy-Schwarz inequality in probability, \eqref{eq:main-step4-1} turns into
  \begin{eqnarray*}
    |[\text{L.~H.~S. of \eqref{pf:main-step5b}}]|&\lesssim&
    \expec{\delta(x{-}\ee_i{-}y\mod
      L)+|\nabla_{x_i}\nabla_y
      G_L(x{-}\ee_i,y)|^4}^{1/4}\\
    &&\qquad\times\,\expec{|[\nabla\phi_j(y)+\ee_j]_y|^4}^{1/4}.
  \end{eqnarray*}
The first term of the r.~h.~s. is estimated by the annealed estimate
  \eqref{eq:annealed2} of Lemma~\ref{L:MO}.
For the second term we appeal
to \eqref{referee} with $q=4$, and to the bound \eqref{eq:unif-bound-corr2}
in Lemma~\ref{L:GNO1-1} on the quartic moment of $\nabla\phi_j$:
\begin{equation}\label{eq:15}
\expec{|[\nabla\phi_j(y)+\ee_j]_y|^4} \,\leq \, 2^4\expec{\big|\nabla\phi_{j}(y)+\ee_j\big|^4 } \,\lesssim \,\expec{|\nabla\phi_{j}|^4 } +1\,\lesssim\, 1.
\end{equation}
The desired estimate \eqref{pf:main-step5b} follows.

\smallskip

\noindent We then turn to \eqref{eq:7}.
Based on \eqref{eq:main-step4-2}, we first bound the
  l.~h.~s. in \eqref{eq:7} by the sum of three terms:
  \begin{align}\label{eq:splitting}
    \begin{split}
      &|[\text{L.~H.~S. of \eqref{eq:7}}]|\,\lesssim\,
      \expec{|\nabla_yG_L(x,y)|^2|F_1(y)|^2}^{1/2}\\
      &\qquad\qquad+\,\expec{|\nabla_yG_L(x,y)|^2
        |F_2(y,x)|^2}^{1/2}\,+\,\expec{|\nabla_yG_L(x,y)|^2 |F_3(y)|^2}^{1/2}.
    \end{split}
  \end{align}
  We estimate the first and third terms on the r.~h.~s. by appealing to the quenched estimate of
  Lemma~\ref{lem:quenched} and Jensen's inequality in probability:
  \begin{eqnarray*}
    \lefteqn{\expec{|\nabla_yG_L(x,y)|^2|F_1(y)|^2}^{1/2}+
      \expec{|\nabla_yG_L(x,y)|^2|F_3(y)|^2}^{1/2}}&&\\
    &\lesssim& (|x-y\mod L|+1)^{2-d-\gamma}\left(\,
      \expec{|\nabla z(y)|^2}^{1/2}+\max_{j=1,\ldots,d}\expec{\phi_j^2}^{{1}/{2}}|\nabla^2 u_0(y)|\right).
  \end{eqnarray*}
  Due to \eqref{eq:unif-bound-corr1} in
  Lemma~\ref{L:GNO1-1} the second moment of $\phi_{j}$ is bounded by $\mu_d(L)$, and thus we can control the r.~h.~s. by
  $I_1$. It remains to estimate the second term in \eqref{eq:splitting}. By appealing to the definition of
  $F_2$, the Cauchy-Schwarz inequality in probability, and \eqref{eq:15} we have
  \begin{equation*}
    \expec{|\nabla_yG_L(x,y)|^2
      |F_2(y,x)|^2}^{1/2}\,\lesssim\,\expec{|\nabla_yG_L(x,y)|^4}^{{1}/{4}}|\nabla
    u_0(x)-\nabla u_0(y)|.
  \end{equation*}
  In view of the annealed estimate \eqref{eq:annealed} of
  Lemma~\ref{L:MO} the r.~h.~s. is controlled by $I_2$ as desired.
 
\medskip
  
\step{5e} Proof of \eqref{pf:main:est3-alt} and thus \eqref{pf:main:est3}.

\medskip

\noindent
By combining \eqref{eq:7} and \eqref{pf:main-step5b} we get
\begin{eqnarray*}
  \lefteqn{[\text{L.~H.~S. of \eqref{pf:main:est3-alt}}]}&&\\
  &\lesssim&
  \left(\sum_{x\in\T_L}\left(\sum_{y\in\T_L}(|y-x\mod L|+1)^{2(1-d)-\gamma}\left(\expec{|\nabla
      z(y)|^2}^{{1}/{2}}+\mu_d(L)|\nabla^2 u_0(y)|\right)\right)^2\right)^{{1}/{2}}\\
  &&+
  \left(\sum_{x\in\T_L}\left(\sum_{y\in\T_L}(|y-x\mod L|+1)^{1-2d}|\nabla u_0(x)-\nabla u_0(y)|\right)^2\right)^{{1}/{2}}.
\end{eqnarray*}
Since $\gamma>0$, the discrete convolution kernel $x\mapsto (|x\mod
L|+1)^{2(1-d)-\gamma}$ has \mbox{$\ell^1(\T_L)$-norm} bounded
independently of $L$ (even for $d=2$), i.~e. \mbox{$\sum_{\T_L} (|x\mod L|+1)^{2(1-d)-\gamma}\lesssim1$}. Hence, by the convolution estimate w.~r.~t. the \mbox{$\ell^2(\T_L)$-norm}, the first term on the r.~h.~s. is controlled by the r.~h.~s.  of \eqref{pf:main:est3-alt}.  It remains to treat the second sum and
suffices to show that
\begin{align}\label{eq:2}
  \begin{split}
    &\left(\sum_{x\in\T_L}\left(\sum_{y\in\T_L}(|y-x\mod L|+1)^{1-2d}|\nabla u_0(x)-\nabla u_0(y)|\right)^2\right)^{{1}/{2}}\\
    &\qquad\qquad\,\lesssim\, \mu_d(L)\left(\sum_{\T_L}|\nabla^2
      u_0|^2\right)^{{1}/{2}}.
  \end{split}
\end{align}
By the definition of the discrete gradient and periodicity we have for $i=1,\ldots,d$ 
\begin{equation*}
  \left(\sum_{x\in\T_L}|\nabla u_0(x+\ee_i)-\nabla u_0(x)|^2\right)^{{1}/{2}}=
    \left(\sum_{\T_L}|\nabla_i\nabla u_0|^2\right)^{{1}/{2}},
\end{equation*}
which combined with the triangle inequality and periodicity yields
\begin{equation}\label{eq:3}
  \left(\sum_{x\in\T_L}|\nabla u_0(x+z)-\nabla u_0(x)|^2\right)^{{1}/{2}}\,\lesssim\,|z\mod L|
  \left(\sum_{\T_L}|\nabla^2 u_0|^2\right)^{{1}/{2}}
\end{equation}
for all $\Z^d$. 
We then use the triangle inequality in the form of
$$
\left( \sum_x \left( \sum_z X_{xz}\right)^2\right)^{1/2}
\,\leq \,  \sum_z \left(\sum_x X_{xz}^2 \right)^{1/2},
$$
so that
\begin{eqnarray*}
  \lefteqn{[\text{L.~H.~S. of \eqref{eq:2}}]}&&\\
  &=&\left(\sum_{x\in\T_L}\left(\sum_{z\in\T_L}(|z\mod L|+1)^{2(1-d)}\frac{|\nabla u_0(x)-\nabla u_0(x+z)|}{(|z\mod L|+1)}\right)^2\right)^{{1}/{2}}\\
  &\stackrel{\text{$\Delta$-inequality}}{\leq}&\sum_{z\in\T_L}\left((|z\mod
    L|+1)^{4(1-d)}\sum_{x\in\T_L}\frac{|\nabla u_0(x)-\nabla u_0(x+z)|^2}{(|z\mod L|+1)^2}\right)^{{1}/{2}}\\
  &\stackrel{\eqref{eq:3}}{\leq}&\sum_{z\in\T_L}(|z\mod L|+1)^{2(1-d)}\,\left(\sum_{\T_L}|\nabla^2 u_0|^2\right)^{{1}/{2}}.
\end{eqnarray*}
Evaluating the sum in $z$ on the r.~h.~s. yields the claimed estimate \eqref{eq:2}, recalling that $\mu_d(L)=\ln L$ for $d=2$ and $\mu_d(L)=1$ for $d>2$.
This proves \eqref{pf:main:est3-alt}, and therefore the desired estimate \eqref{pf:main:est3}
by Step~5a.

\qed


\subsection{Proof of Corollary~\ref{cor:1}}
Estimate \eqref{eq:coro1} is a direct consequence of \eqref{eq:main} and the estimate
\begin{equation*}
  \max_{j=1,\ldots,d}\expec{|\phi_j|^2}\,\lesssim \, \left\{
\begin{array}{ll}
\ln L& \mbox{ for }d=2\\
1&\mbox{ for }d>2
\end{array}
\right.
\end{equation*}
of Lemma~\ref{L:GNO1-1}.
\qed
\appendix

\section{Proofs of the auxiliary lemmas}\label{sec:proof2}

\subsection{Proof of Lemma~\ref{L1}}

We adapt the arguments of \cite[Proof of Lemma~7]{Gloria-Neukamm-Otto-11a}, as for the proof of \cite[Lemma~3]{Gloria-Otto-09b} starting from  \cite[Lemma~2.3]{Gloria-Otto-09}.
We first introduce a couple of notations: Let
$\{y_n\}_{n=1,\ldots,N}$, $N:=L^d$, be an
enumeration of 
$\T_L$. For $n = 1,\cdots,N$ define the average
\begin{equation*}
  \expec{\cdot}_{\leq
  n}:=\prod_{1\leq k\leq n}\int_{\Omega_0}\beta(d\aa(y_k)).
\end{equation*}
Set $\zeta_n:=\expec{\zeta}_{\leq n},\,\tilde\zeta_n:=\expec{\tilde \zeta}_{\leq n}$ for $n\geq 1$ and
$\zeta_0:=\expec{\zeta}_{\leq 0}:=\zeta,\,\tilde \zeta_0:=\expec{\tilde \zeta}_{\leq 0}:=\tilde \zeta$. 
W.~l.~o.~g. we assume that both $\zeta$ and $\tilde \zeta$ have zero expectation.
We split the proof into two steps.
  
\medskip

\step 1 Martingale decomposition.
\medskip
  
\noindent We claim that
\begin{equation}
\label{eq:27}
\expec{\zeta\tilde \zeta}=\sum_{n=1}^{N} \expec{(\zeta_n
  -\zeta_{n-1})(\tilde \zeta_n-\tilde \zeta_{n-1})}.
\end{equation}
Here comes the argument: Since $\{y_1,\ldots,y_N\}=\T_L$ we
have  $\expec{\cdot}=\expec{\cdot}_{\leq N}$ by \eqref{eq:def-iidens},
and thus $\zeta_N=\expec{\zeta}=0$, $\tilde \zeta_N=\expec{\tilde
  \zeta}=0$ and $\zeta_N\tilde \zeta_N=0$. Moreover, by construction
we have $\zeta_0=\zeta$ and $\tilde
\zeta_0=\tilde \zeta$, and thus 
\begin{equation}
\label{eq:26}
\zeta\tilde \zeta=\sum_{n=1}^N\zeta_{n-1}\tilde \zeta_{n-1}-\zeta_n\tilde \zeta_n.
\end{equation}
The identity \eqref{eq:27} then follows from taking the expectation of \eqref{eq:26} provided we prove that
\begin{equation}\label{eq:app-felix}
\expec{(\zeta_n -\zeta_{n-1})(\tilde \zeta_n-\tilde \zeta_{n-1})}\,=\,\expec{\zeta_{n-1}\tilde \zeta_{n-1}}-\expec{\zeta_n\tilde \zeta_n}.
\end{equation}
This last identity can be seen as follows: By definition we have
\begin{equation*}
  \tilde\zeta_n=\expec{\tilde\zeta}_{\leq n}=\prod_{1\leq k\leq n}\int_{\Omega_0}\tilde\zeta\,\beta(d\aa(y_k))=\int_{\Omega_0}\tilde\zeta_{n-1}\,\beta(d\aa(y_n)).
\end{equation*}
Since $\zeta_n$ does not depend on $y_1,\ldots, y_n$, we have
\begin{equation*}
  \zeta_n\tilde\zeta_n=\int_{\Omega_0}  \zeta_n\tilde\zeta_{n-1}\,\beta(d\aa(y_n)).
\end{equation*}
Integrating both sides w.~r.~t. $\aa(y_k)$, $k\neq n$,
yields $\expec{\zeta_n\tilde\zeta_n}=\expec{\zeta_n\tilde\zeta_{n-1}}$ and thus by symmetry
$\expec{\zeta_n\tilde\zeta_n}=\expec{\zeta_{n-1}\tilde\zeta_{n}}$, 
so that we obtain \eqref{eq:app-felix}.

\medskip

\step 2 Conclusion.
\medskip
   
\noindent From \cite[Step~2, Proof of
Lemma~7]{Gloria-Neukamm-Otto-11a}, since $\expec{\cdot}$ is an
$L$-periodic i.~i.~d. measure, for all $n\in\{1,\ldots,N\}$ we have
\begin{equation*}
\expec{(\zeta_{n-1}-\zeta_n)^2}\leq
\Expec{\left(\frac{\partial\zeta}{\partial y_n}\right)^2}, \quad \expec{(\tilde \zeta_{n-1}-\tilde \zeta_n)^2}\leq
\Expec{\left(\frac{\partial\tilde \zeta}{\partial y_n}\right)^2}.
\end{equation*}
Hence the claim follows from \eqref{eq:27} and Cauchy-Schwarz' inequality.

\qed

\subsection{Proof of Lemma~\ref{lem:quenched}}

It suffices to consider the case $y=0$, since
$G_L(x,y;\aa)=G_L(x-y,0;\aa(\cdot-y))$ and due to the fact that the
asserted estimate depends on $\aa$ only through its ellipticity constant $\lambda$. We shall write $G_L(x):=G_L(x,0;\aa)$ for brevity.
We divide the proof into three steps.
First we derive quenched estimates on the Green's function for $d>2$ from the quenched estimates  of \cite[Lemma~24 and (219)]{Gloria-Neukamm-Otto-11a} on the corresponding parabolic Green's functions.
In the second step we prove the desired estimate for $d>2$ using the De Giorgi-Nash-Moser H\"older estimate, 
see \cite[Theorem~5.2]{Mosconi-01} in the discrete setting.
We then deduce the estimate for $d=2$ from the estimate for $d=3$ using an argument by Avellaneda and Lin, see
\cite[Proof of Theorem~13]{Avellaneda-Lin-87}. 

\medskip

\step 1 Quenched estimate on $G_L$ for $d>2$. 

\medskip

\noindent 
Let $G_L'(t,y)$ denote the periodic, parabolic Green's
function considered in \cite[Definition~2]{Gloria-Neukamm-Otto-11a};
it is characterized as follows: For all
$x\in\Z^d$ the function $G_L'$ is the unique $C^\infty(\R,\ell^\infty(\Z^d))$ solution of
\begin{equation}\label{eq:quenched-step1-parabolic}
\left\{
\begin{array}{rcl}
\partial_t G_L'(t,x)-\nabla \cdot \aa \nabla G_L'(t,x)&=&0,\\
G_L'(0,x)&=& \delta(x \mod L).
\end{array}
\right.
\end{equation}
Since both the initial data and the coefficient field are
$L$-periodic, $G_L'$ is in fact of class $C^\infty(\R,\ell^1(\T_L))$.
The elliptic Green's function $G_L$ can be formally obtained by integrating in time the mean free version of the parabolic Green's function:
\begin{equation*}
G_L(x)\,=\,\int_0^\infty (G_L'(t,x)-L^{-d})dt.
\end{equation*}
To turn this into a rigorous argument one needs to prove that $t\mapsto G_L'(t,x)-L^{-d} \in L^1(0,+\infty)$.
This is a direct consequence of \cite[eq. (219) in the proof of the discrete
version of Lemma~24]{Gloria-Neukamm-Otto-11a}: For all $\alpha>0$ we have
\begin{eqnarray*}
0\leq   G_L'(t,x)&\lesssim& (t+1)^{-\frac{d}{2}}\Big(\frac{(|x \mod L|+1)^2}{t+1}+1\Big)^{-\frac{\alpha}{2}}\quad\mbox{for}\;t\lesssim L^2,\label{L2.11}\\
  |G_L'(t,x)-L^{-d}|&\lesssim& L^{-d}\exp(-c_0\frac{t}{L^2})\quad\mbox{for}\;t\gtrsim L^2.\label{L2.10}
\end{eqnarray*}
For $d>2$, this directly yields 
\begin{equation}\label{eq:pr-quenched-2}
|G_L'(x)|\,\leq \,\int_0^\infty |G_L(t,x)-L^{-d}|dt \, \lesssim \, (1+|x\mod L|)^{2-d}.
\end{equation}

\medskip

\step 2 Quenched estimate on $\nabla G_L$ for $d>2$.
\medskip

\noindent We now deduce \eqref{eq:quenched} from
\eqref{eq:pr-quenched-2} for $d>2$. Recall that it suffices to prove
the estimate for $y=0$ and that we write $G_L(x)=G_L(x,0)$
for brevity. Fix a radius $2\leq R \sim 1\ll L$.
For $|x\mod L|\leq
2R$ the r.~h.~s. of \eqref{eq:quenched} is of order $1$ so that the estimate
directly follows from the combination of \eqref{eq:pr-quenched-2}
with the discrete estimate $|\nabla_i G_L(x)|\leq |G_L(x+e_i)|+|G_L(x)|$. Hence, it suffices to consider the case
$|x\mod L|\geq 2R$. 
For all $y'\in \T_L$ and $r>0$, set $B_r(y')=\{y''\in \T_L, |y'-y''\mod L|\leq r\}$.
Using the elementary inequality
$$
|\nabla G_L(x)|\,\leq \, \sqrt{d}\osc{B_1(x)}{G_L},
$$
this turns into an estimate of the oscillation of $G_L$ on $B_1(x)$.
Noting that $G_L$ satisfies  
\begin{equation}\label{eq:quenched-step2-elliptic}
  \nabla^*\aa\nabla G_L=L^{-d}
\end{equation}
in the set $B_{|x \mod L|/2}(x)$ (which does not contain $0 \mod L$),
one may appeal to the De Giorgi-Nash-Moser H\"older estimate. In particular, 
by \cite[Theorem~5.2]{Mosconi-01}, there exists $0<\gamma<1$ depending only on $\lambda$ and $d$ such that
\begin{equation}\label{harnack-1}
  \osc{B_{1}(x)}{G_L}\, \lesssim \, \left(|x\mod
  L|^{-\gamma}\max_{B_{|x \mod L|/2}(x)} |G_L|\right)+ L^{-d},
\end{equation}
Combined with \eqref{eq:pr-quenched-2}, this yields \eqref{eq:quenched} for $d>2$.

\medskip

\step 3 Quenched estimate on $\nabla G_L$ for $d=2$.
\medskip

\noindent 
Following \cite[Proof of Theorem~13]{Avellaneda-Lin-87} we derive the
estimates for $d=2$ from the estimate for $d=3$. To distinguish
quantities in different dimensions we use the superscripts $^{(2)}$
and  $^{(3)}$; e.g. $\T_L^{(2)}$ denotes the  $2$-dimensional
torus. To a given two-dimensional, $L$-periodic
coefficient field $\aa^{(2)}\in\Omega_L^{(2)}$ we associate a three-dimensional coefficient-field
$\aa^{(3)}\in\Omega_L^{(3)}$ via
$$
\aa^{(3)}(x,x_3)\,:=\,\dig{[\aa^{(2)}(x)]_1, [\aa^{(2)}(x)]_2,1}\qquad(x,x_3)\in\Z^2\times\Z=\Z^3.
$$
In the following we use the shorthand notation
$G_L^{(2)}(x):=G_L^{(2)}(x,0;\aa^{(2)})$ and $G_L^{(3)}(x,x_3):=G_L^{(3)}((x,x_3),(0,0);\aa^{(3)})$.
It is elementary to check that
$$
G_L^{(2)}(x)\,=\,\sum_{x_3\in([0,L)\cap\Z)} G_L^{(3)}(x,x_3),
$$
and thus
$$
\nabla_iG_L^{(2)}(x)\,=\,\sum_{x_3\in([0,L)\cap\Z)}
\nabla_iG_L^{(3)}(x,x_3)\qquad\text{for }i=1,2.
$$
The quenched estimate \eqref{eq:quenched} for $d=3$ then turns into 
\begin{eqnarray*}
|\nabla_iG_L^{(2)}(x)|&\leq
&\sum_{x_3\in([0,L)\cap\Z)}|\nabla_i G_L^{(3)}(x,x_3)| \\
&\lesssim & \sum_{x_3\in([0,L)\cap\Z)}(1+|(x,x_3) \mod L|)^{-1-\gamma}\\
&\lesssim & \sum_{x_3\in([0,\infty)\cap\Z)} (1+|x \mod L|+x_3)^{-1-\gamma} \\
&\lesssim  & (1+|x\mod L|)^{-\gamma}\lesssim (1+|x\mod L|)^{-\gamma}
\end{eqnarray*}
which is nothing but \eqref{eq:quenched} for $d=2$.
\qed

\appendix

\end{document}